# UNDERSTANDING TEACHERS' PROFESSIONAL DEVELOPMENT THROUGH THEIR INTERACTIONS WITH RESOURCES: A MULTILINGUAL PROJECT

Luc Trouche

French Institute of Education, Ecole normale supérieure de Lyon, France

Luc.Trouche@ens-lyon.fr

*The Documentational Approach To Didactics (DAD) aims to study teachers' professional development through their interactions with their resources for/from teaching. It has been introduced in the French community of didactics of mathematics in 2007, then extended at an international level. It has been introduced as an entry in the Encyclopedia of Mathematics Education in 2020. The DAD-Multilingual project (2020-2021) is dedicated to gather and confront the translations of this entry towards 14 languages. The project main goals are: making available – for students as well as for researchers - a presentation of DAD in various languages; deepening the DAD concepts themselves in thinking their possible instantiations in different languages; questioning the translation processes; and questioning the notion of resource itself, resource for/from teaching. The lecture presents this project, and draws some lessons from its first steps.*

In this lecture, we want to present an on-going project dedicated to better understand mathematics teachers' professional development through the lens of their interactions with a diversity of resources. In the first part, we will introduce the so-called Documentational Approach to didactics (DAD), which has been developed for about 10 years. In a second part, we will present the DAD-Multilingual project, aiming to deepen this approach through its adaptation towards different social, curricular and linguistic contexts. In the third part, we will present the feedback of the scientific committee of this project, allowing to better situating the scope of the project and the ways for its development. In the fourth part, we will present the preliminary results, and will conclude in drawing some perspectives.

## 1. DAD, towards a 'resource' approach to mathematics education.

The documentational approach to didactics (DAD) has been introduced by Ghislaine Gueudet and Luc Trouche (Gueudet & Trouche, 2009), and has been developed further in joint work with Birgit Pepin (Gueudet, Pepin & Trouche, 2012). We will just introduce here the main concepts of this approach; more information may be found in the DAD entry (Trouche, Gueudet & Pepin, 2020) of the second edition of the Encyclopedia of Mathematics Education edited by Stephen Lerman.

DAD is originally steeped in the French didactics tradition in mathematics education (Artigue et al., 2019), where concepts such as *didactical situation*, *institutional constraint* and *scheme* are central. At the same time it also leans on socio-cultural theory, including notions such as *mediation* (Vygotsky, 1978) as constitutive of each cognitive process. Moreover, the approach has also been developed due to the emerging digitalization of information and communication, which asks for new conceptualizations. This digitalization and the development of Internet had indeed strong consequences: ease of quick access to many resources and of communication with many people. This necessitated a complete metamorphosis of thinking and acting, particularly in education: new balances between *static* and *dynamic* resources, between *using* and *designing* resources, between *individual* and *collective* work (Pepin, Choppin, Ruthven, & Sinclair, 2017). Taking into account these phenomena, DAD proposed a change of paradigm by analyzing teachers' work through the lens of "resources" for and in teaching: what they prepare for supporting their classroom practices, and what is continuously renewed by/in these practices. This sensitivity to *resources* meets Adler's (2000) proposition of "think[ing] of a resource as the verb re-source, to source again or differently" (p. 207). Retaining this point of view, DAD takes into consideration a wide spectrum of resources that have the potential to resource teacher activity (e.g. textbooks, digital resources, email exchanges with colleagues, or student worksheets), resources *speaking to the teacher* (Remillard, 2005) and supporting her/his engagement in teaching.

During the interaction with a particular resource, or sets of resources, teachers develop their particular *schemes of usage* of these resources. The concept of "scheme" (Vergnaud, 1998) is central in DAD. It is closely linked with the concept of "class of situations", which are, in our context, a set of professional situations corresponding to the same aim of the activity (for example, introducing a given mathematical property for a given grade). For a given class of situations, a teacher develops a more or less stable organization of his/her activity, that is a scheme. A scheme has four components:

- The aim of the activity;
- Rules of action, of retrieving information and of control;
- Operational invariants, which are elements, often implicit, of knowledge guiding the activity;
- Possibilities of inferences, meaning of adaptation to the variety of situations.

Over the course of his/her activity, a teacher enriches his/her schemes, e.g., integrating new rules of actions, or s/he can develop new schemes. Schemes are likely to be different for different teachers, although they may use the same resources, depending on their dispositions and prior knowledge.

The resources and the scheme, developed by a given teacher for facing a given class of situations, make up a document. The process of developing a document has been coined *documentational genesis* (Figure 1). The 'use' of resources is an interactive and potentially transformative process. This process works both ways: the affordances of the resource/s influence teachers' practice (that is the *instrumentation* process), as the teachers' dispositions and knowledge guide the choices and transformation processes between different resources (that is the *instrumentalisation* process). Hence, the DAD emphasizes the dialectic nature of the teacher-resource interactions combining instrumentation and instrumentalisation. These processes include the design, re-design, or 'design-in-use' practices (where teachers change a document 'in the moment' and according to their instructional needs).

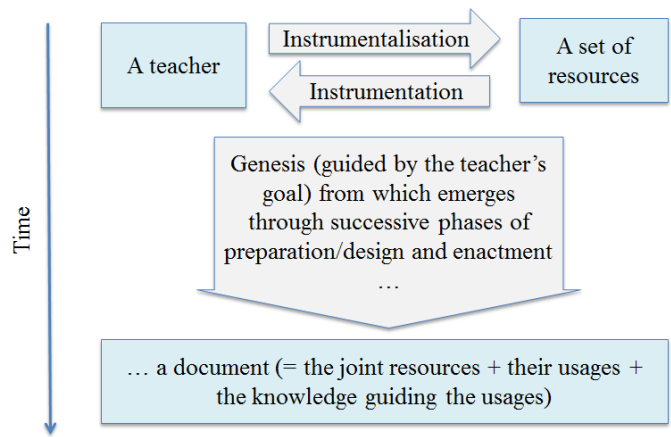

**Figure 1. A representation of a documentational genesis**

The set formed by all the resources used by the teacher is named his/her *resource system*. These resources are associated with schemes of usage, forming documents. The documents developed by a teacher also form a system, called the *document system* of the teacher. Its structure follows the structure of the class of situations composing the professional activity of the teacher. When teachers share their documentation work, for example in a group preparing lessons collectively, they may also develop a shared resource system (Gueudet, Pepin & Trouche, 2012). Nevertheless, the different members of the group can develop different schemes for the same resource, resulting in different documents.

We have then presented the main concepts grounding DAD. Since its introduction, this approach has been used in a variety of contexts, in Ph.D. and research projects. The Re(s)source international conference, held in 2018 in Lyon (https://resources-2018.sciencesconf.org/), gathering 130 people from 30 countries, gave a good image of the extension of the French original field. This cultural and linguistic diversity was understood as a potential richness for deepening the concepts at stake:

- One of the sessions of the young researchers workshop, held during this conference, was dedicated to « Naming systems[1] used by secondary school teachers to describe their resources and their documentation work », meaning the structured set of words used by teachers, in their own language, for describing their resource systems;
- And, in my final conference (Trouche, 2019), among the 10 research programs that I proposed for developing DAD, two of them addressed linguistic issues: the first one,

---

[1] The notion of "naming system" was inspired by the Lexicon project. "[This project involved] nine countries (Australia, Chile, China, Czech Republic, Finland, France, Germany, Japan, USA), and in each country a team of mathematics education researchers and experienced mathematics teachers. In this project, we consider that our experiences of the world and reflection on those experiences are mediated and shaped by available language, and that the use of English as *lingua franca* for international communication substantially limits what can be expressed and shared. The goal of the project is thus to document and compare the naming systems employed in mathematics teacher communities in the nine countries to describe the objects and events in their classrooms, in order to expand on the variety of constructs available for the purpose of theorizing about classroom practice and for identifying the characteristics of accomplished practice" (Artigue et al., 2019). But while the Lexicon project looked at the naming systems used to describe classroom activities, the young researchers workshop looked at the naming systems used to describe teachers' interactions with resources, before, during, and after class.

"Conceiving a DAD living multi-language glossary", and the last one "Contrasting naming systems used by teachers in describing their resources and documentation work, towards a deeper analysis of teachers' resource systems."

These reflections, among others, lead to the DAD-Multilingual project.

## 2. The DAD-Multilingual project, deepening a theoretical approach through its adaptation to a diversity of contexts

We describe here the origin of the project, the actors involved, and the translating processes. Any project is actually born from the convergence of a set of phenomena; and responds to a set of needs. It is indeed the case for the DAD-Multilingual project, being:

- The result of my personal experience, as a French native speaker having to go back and forth between English and French: introducing first DAD in French (Gueudet & Trouche, 2008), then in English (Gueudet & Trouche 2009); writing in English the entry DAD for the Encyclopedia of Mathematics Education (Trouche, Gueudet, & Pepin, 2020), then translating it in French (Trouche, Gueudet, Pepin, & Aldon, 2020). Doing so, I had in mind the Tuareg proverb: "travelling is going from oneself to oneself through others"…
- The result of my PhD supervisions experience: I realize, for example, that students from some countries had never the occasion to express themselves, in the frame of their studies, in their own language: "In spite of the numerous calls from education and language specialists, many countries still use the languages of wider communication instead of their native languages […] As a result, students are often required to learn subject material in the language of a former power; a language in which they may not be proficient "(Quigley et al., 2011). The DAD-Multilingual project appears then as a necessity for addressing equities issues;
- The result of interacting with researchers in various contexts (Algeria, Brazil, China, Japan, Lebanon, Mexico, Netherlands, Norway, or Senegal): these interactions have evidenced the fact that dealing with DAD in each new context (theoretical, cultural, curricular a well a linguistic) leads to new questions and potential enrichment. As it was said for didactics (Arcavi et al., 2016), DAD "goes travelling"…
- The result of developing an approach grounded on teachers' work with resources, involving naturally a diversity of supports and languages;
- The result of working over a long period with Ghislaine Gueudet and Birgit Pepin, committed both in international projects, Birgit having herself a long experience of crossing linguistic boarders…

These interactions evidenced also the need for enlightening the complex metaphoric structure developed by DAD, that could appear as a characteristic of the French community of mathematics education: "Despite the broad dispersion and wide-ranging accomplishments of didactique over the past decades, it has not had the influence outside the Francophone world that one might have expected […] part of the communication problem is that didactique carries some heavy baggage stemming largely from the language it employs and its metaphors in particular […] Didactique, in creating a precise vocabulary for its work, has made extensive use of the fundamental metaphoric structure identified by Pimm [1988, 2010], generating terms that need careful exegesis before they are used. Anglophones may find that English versions of those terms come laden with extra baggage that makes them difficult to interpret correctly." (Jimmy Kilpatrick in Arcavi et al., 2016)[2].

Finally, the presentation of DAD in the Encyclopedia of Mathematics Education gave us (Ghislaine, Birgit and me) the opportunity of a conceptual reversal. Each Encyclopedia, since Diderot and d'Alembert's work (1760), rests on a fundamental objective: making available all the knowledge of the world in a given place (a series of books) and a given language. Our project, reversing this objective, was to make available a small piece of knowledge in a diversity of languages, with the idea that this diversity will contribute to better understand the piece of knowledge at stake. Thus was born the DAD-Multilingual project (https://hal.inria.fr/DAD-MULTILINGUAL, aimed at adapting the Documentational Approach to Didactics entry into a diversity of languages.

The translating process involved 14 languages (in addition to English), actually the languages represented in the Re(s)source 2018 international conference (Gitirana *et al.* 2018): Arabic, Chinese, French, German, Greek, Hebrew, Hungarian, Italian, Japanese, Norwegian, Portuguese, Spanish, Turkish and Ukrainian. These 14 languages offer both elements of proximity (as for the roman languages: French, Italian, Portuguese and Spanish) and elements of distance, for example between European languages and Chinese one, leading to conceptualize differences of languages and of thought (Jullien, 2015).

The goals of the project, as announced on its website (https://hal.inria.fr/DAD-MULTILINGUAL), are the following ones:

- Making available a presentation of DAD in various languages, allowing the students and the researchers interested to refer to it in their own language;
- Deepening the DAD concepts themselves in thinking their possible instantiations in different languages;
- Questioning the translating process itself;
- Beyond the frame of DAD, questioning the notion of 'resource', resource for/from teaching;
- Designing (in a later step) a multilingual glossary of DAD.

The project involves a set of translators, reviewers (at least one translator and one reviewer for each language) and a scientific committee (see § 3). It would be excessive to say that they formed a community at the start of the project. Actually, there is not a DAD community, no regular event or specific journal allowing labelling such a scientific group. Since its beginning, in 2009, DAD develops as an « approach », with blurred boundaries, acting as a theoretical workshop for studying teachers-resources interactions, complementing or questioning already well established theoretical frameworks. During these 10 years, DAD has attracted PhD students, and researchers, in the frame of projects around mathematics teaching resources: e.g., in France, ReVEA (https://www.anr-revea.fr); in Europe, MC2 (http://www.mc2-project.eu); internationally, the French-Chinese joint project MaTRiTT (http://ife.ens-lyon.fr/ife/recherche/groupes-de-travail/matritt-joriss). The translators and reviewers have a diversity of links to DAD: interested as prospective users, or effective users, or co-designers (in particular, PhD students have enriched DAD with new concepts in their theses). They all are native speakers for the targeted language of a given translation; and they are sometimes go-between different languages, for historical reasons (e.g., Arabic-French in Lebanon) or PhD reasons (e.g. Chinese students having done their PhD in a frame of a co-supervision, using French and Chinese for collecting data; and English for writing their thesis), or a mix of these reasons (see Window 1).

---
[2] Thanks to Tommy Dreyfus who, after his reading of a preliminary version of this paper, draws my attention on this Kilpatrick's contribution.

---

**Window 1 - The translator-reviewer pair in the Spanish language case**

The translator was Ulises Salinas-Hernández, and the reviewer Ana Isabel Sacristán. Both have been members of the Department of Mathematics Education in Cinvestav-IPN (Mexico):

- Ulises obtained his PhD from Cinvestav; then has been doing a two-year post-doctorate at the ENS de Lyon with Luc Trouche, reflecting on theoretical networking, crossing DAD and the semiotic approach (Radford, 2008). His stays at ENS de Lyon gave him the opportunities to contrast the naming systems used by Mexican and Chinese teachers (Wang, Salinas & Trouche, 2019);
- Ana Isabel did her PhD at the University of London with Richard Noss and partially with Celia Hoyles. She has a long history of interacting with Luc Trouche, first at the ENS de Lyon during a three-month scientific stay in 2012, as well as in Cinvestav during a two-month scientific stay of Luc Trouche in 2017; these interactions gave rise to several papers (e.g., Trouche, Drijvers, Gueudet & Sacristan, 2013).

These close interactions have allowed flexible discussions on the translating process. This is a specific case for the translator-reviewer pair, and other cases can be found in the project, more or less close to DAD history.

---

The method of the translating process was as follows:

- English was the interface language (of course, this could constitute a bias for the on-going discussions);
- From the beginning, it was clear that the objective was not to produce a translation that was as close as possible to the original text, but instead to: make DAD understandable in a specific cultural, curricular and linguistic context (that means bridging it, if possible, with other frames well known of the targeted audience); enrich DAD in questioning its concepts when translating it;
- Issues that arose in the translating process were shared by each translator-reviewer pair, who had to fill, in English, a 'Translating issues report'. For such a report, a model was proposed (see Window 2), but such a model could be adapted according to the needs of the translator-reviewer pair.

---

**Window 2 – The model (to be adapted) for the Translating issues report**

Language:    Translator:    Reviewer:
Sources: English version and other linguistic versions?

- In a few lines, could you describe the main issues that emerged when translating the DAD entry or when interacting with the reviewer? Issues linked to the context (social, cultural, or curricular); issues linked to the concepts at stake; issues linked to the

> - vocabulary
> - Certain concepts raised difficulties, or discussions between the translator and the reviewer. We suggest that you explain these difficulties, and the choices you have made, for the notions of resource, document and for about three other notions, which seemed more particularly complex: Possible translations, and associated definitions (in English) - Final choice, and motivation - Scientific references using this word in the targeted language
> - Other issues that you would like to share

Each translation was considered as an element of a collection, integrated in a French scientific Open Archive website (see Window 3).

> **Window 3 – The presentation of the Chinese translation on the Open Archive Website**
>
> **文献纪录教学论**
>
> Luc Trouche, Ghislaine Gueudet, & Birgit Pepin1
> 中文翻译：王重洋（Chongyang Wang）
> 中文审阅：徐斌艳（Binyan Xu）
>
> **摘要**
> 文献纪录教学论是数学教育百科全书的入门词条(Trouche, Gueudet & Pepin 2018)。词条于2020年更新 (Trouche, Gueudet & Pepin 2020)。本文是中文译本，是14种语言译本之一 (https://hal.archives-ouvertes.fr/DADMULTILINGUAL）。
>
> 文献纪录教学论是数学教育领域的一种教学理论。该理论创立的初衷是通过研究教师和资源之间的互动（包括使用和设计）来理解教师专业发展。本文主要阐述该理论的起源，理论背景，核心概念和相关的方法论。为达到理据结合的效果，我们将结合不同的研究项目作为案例诠释上述内容。本文面向的读者群是研究者，以及对入门文献纪录教学论有兴趣的非专家型读者（比如硕士研究生）。
>
> **Abstract**
> The 'Documentation Approach to Didactics' is an entry of the Encyclopedia of Mathematics Education (Trouche, Gueudet & Pepin 2018). This entry has been updated in 2020 (Trouche, Gueudet & Pepin 2020). This article is a Chinese adaptation of this updated version. It is part of a collection, gathering such adaptations in 14 languages (https://hal.archives-ouvertes.fr/DADMULTILINGUAL).
> The documentational approach to didactics is a theory in mathematics education. Its first aim is to understand teachers' professional development by studying their interactions with the resources they use and design in/for their teaching. In this text we briefly describe the emergence of the approach, its theoretical sources, its main concepts and the associated methodology. We illustrate these aspects with examples from different research projects. This synthetic presentation is written for researchers, but also for non-specialists (e.g. master students) interested in a first discovery of the documentational approach.
>
> **关键词**
> 课程材料；电子资源；文献起源；操作不变量（行动知识）；资源系统；教学资源；教师集体工作；教师专业发展。
>
> **Keywords**
> Curriculum materials; Digital resources; Documentational geneses; Operational Invariants; Resource systems; Resources for teaching; Teachers' collective work; Teacher professional development.

This website (https://hal.inria.fr/DAD-MULTILINGUAL) gives access to the presentation of the project and its actors; to the set of translations; to the set of Translating issues reports; and to different resources aimed to support the translating processes (most of them coming from the scientific committee); and to the analyses produced over the project (as this current lecture!).

## 3. What we have learnt until now from the feedback of the scientific committee

The scientific committee was at the beginning composed of 5 persons: Jill Adler, Nicolas Balacheff, Rongjin Huang, Janine Remillard and Kenneth Ruthven. They were called upon due to their knowledge in, and interest of: the international community of mathematics education; the resource approach to mathematic education; the semantic issues at stake in each translating process; or/and the interactions between different cultures and languages. From the beginning of the project, they were asked to comment on the way the project was organized (e.g., the model of 'Translating issues report', improved thanks to their comments), and to propose references that could support the reflections on the translating processes. Their full comments can be found on the project website (https://hal.inria.fr/DAD-MULTILINGUAL/page/translation-issues), where we underline what appear as their main contributions.

Kenneth Ruthven draws attention to the source language, with four fundamental questions:

- Why not retain key terms from the source language?
- Why not 'mark' key terms in some way to indicate the specialized usage intended (e.g., reSource)?
- Why not 'mark' key terms in some way to clarify the metaphor (e.g., resource-scheme-document, abbreviated, say, to res-sch-doc)?
- Does a concept (as 'resource system' see Ruthven 2019) need a sharper definition before it can become a key term of DAD? And would that sharper definition point to a more precise term (or phrase)?

In such cases, the support of dictionaries (e.g., the Cambridge Dictionary https://dictionary.cambridge.org/dictionary/ and the already accepted translated terms in the specialized domain (e.g., such as the case of 'scheme' in the field of psychology) should be followed.

Nicolas Balacheff recalls that "The issue of language is not just a question of words, as is too often stated, but of expression and the circulation of meaning" (Balacheff, 2018). The minimal condition for doing this work should be to complete the choice of translated terms with authority quotes attesting to their use, and allowing taking into account the "finesse" of the concepts.

Jill Adler considers that the focus on DAD concepts as isolated words is too narrow. She suggests to take into account the context in which these concepts are used (see also Arcavi et al., 2016; Pepin, 2002; Setati, 2003), and raised the issue of the link between teachers' discourse-resources (Adler, 2012), and the theoretical discourse analyzing them.

Janine Remillard, building on her experience in the Math3Cs project (Remillard, 2019), evokes Osborn's (2004) discussion of different types of equivalences in cross-cultural research and Clarke (2013)'s notions of validity when doing cross-cultural research. Like Jill Adler, she underlines the importance of the context, pointing to how the words themselves are "the tip of the iceberg". For facing these issues, the translation team needs to develop what Andrews (2007) calls *prerequisite intersubjectivity*, leading to a shared understanding of the core concepts (see also Pepin et al., 2019). In this perspective, the design of a multilingual glossary of key terms seems crucial – this is actually one of the objectives of the DAD-Multilingual project, already evoked at the Re(s)source 2018 conference (Trouche, 2019).

For taking into account the link between the cultural context, teachers' words, and the conceptualization of their interaction with resources, the coordinators of the project invited Michèle Artigue, involved in the Lexicon Project (see footnote 1), to join the scientific committee, and so it became composed of 6 persons. Until now, this scientific committee is only composed of members of the mathematics education community: other scientific fields could be, of course, considered (e.g., linguistics, computer science, anthropology, cultural studies…); perhaps to be discussed in a later stage of the project?

## 4. Some preliminary results

As I write the presentation of this lecture (August 20[th] 2020), the project is still on-going. 11 translations (of over 14) have been completed, and the three remaining translations will be achieved before the end of September. This is the first *productive* result of the project, corresponding to its first aim.

Regarding its *constructive* results – what do we learn from the translating process itself – the Translating Issues Reports, still under progress, as well as the interactions within the project, allow drawing some preliminary lessons.

First of all, each translation gave rise to very active processes, often mobilizing several sources, and more actors than the translator and the reviewer. For example, for the Turkish translation, the English and French versions were used, and the translation process was an opportunity to introduce, and to discuss the approach with doctoral students:

> About the translation, we checked the translation together but it seems it is not possible to finish it at the end of the April. Because the sentence type and the explanations are very different from English and French. Doctoral students also find the translation problematic and we are revising it according to their feedback (email from Burcu Nur Basturk, on April 13th 2020)

Second, each translation appeared as a complex process, involving several levels: vocabulary, scientific expressions, and structure of the sentences, as detailed in the Japanese report (Window 4).

---

**Window 4 – Extract of the Japanese Translating issues report**
Takeshi Miyakawa and Yusuke Shinno

After reading all through the translated text of the DAD entry, we found that the text was not really the one we usually write by ourselves in Japanese. One may find that this is the translation, not the original text. This would be due to the difficulties of translation at different levels.

First, at the level of vocabulary, there are many technical terms, which are not used in the ordinary language. We had to create an appropriate Japanese term for the English or French term. This difficulty is not only for the technical terms used in DAD, but also those used in the mathematics education, in the scientific papers in general, or in the ordinary language. For the technical term, we used sometimes the English phonetic expression, and other times the Japanese translated terms. The most difficult term we discussed a lot was the name of approach "Documentational approach to didactics". Even the usual term "approach" was not easy for us to translate.

The use of technical terms is also related to the context of scientific research. In the research on mathematics education in Japan, the scholars often try to use the terms which are comprehensible to others and actually use much less number of technical terms than in the didactics of mathematics in France. Japanese scholars therefore may be surprised with the use of technical terms in this text and sometimes might be uncomfortable with it.

There were also many difficulties of translation at the level of sentence. In Japanese language, the order of terms in a sentence is very different from French and English: for example, the verb is given at the end of the sentence; the subject is not sometimes given; and so forth. Due to this, we had to often split a sentence into several sentences. Further, we consider that the context in which the original English text was written would be a factor that makes our translation alien from the Japanese ordinary text. Some English sentences, which seem self-explanatory would not clearly explain the claim, and Japanese readers may feel the lack of sentences that complement the claim, since they are in the other context.

Third, translating DAD needed to think more globally on the theoretical background of this approach, and on the existing, or not, bridges towards the targeted language, requiring, sometimes, to go through a third language, such as Russian in the case of Ukrainian. In this perspective, the objects and methodologies have to be questioned (see Window 5).

**Window 5 – Extract of the Ukrainian Translating issues report**
Maryna Rafalska and Tetyana Pidhorna

There are no translations of the works of G. Brousseau, G. Vergnaud, and Y. Chevallard in Ukrainian. Thus, we faced the difficulties in translating the main notions used in their theories (e.g. "milieu", "scheme", "operational invariants", "savoir à enseigner" and "savoir enseigné", etc.). For the translation in Ukrainian of the notion of scheme we used Ukrainian articles [11, 12] that refer to the work of J. Piaget and the Russian translations of J. Piaget works [10]. We also used the Russian translation [12] of Rabardel's instrumental approach for translating the notions of instrumentation, instrumentalisation, and instrumental genesis. The Russian terms were translated in Ukrainian using the dictionaries that provide the meaning of these words and then we found the Ukrainian analogues (also using the dictionary to confirm that the meanings of found words are the same).

Translating the DAD showed us the differences in research objects and methodologies in different cultural contexts, e.g. French didactics of mathematics and Ukrainian method of teaching and learning mathematics (metodika navchannya matematiky). Thus, we noticed that in Ukrainian method of teaching and learning mathematics much less attention, compared to French didactics of mathematics, is given to the study of psychological constructs (e.g. scheme) that influence on teachers' choices as well as to transposition of mathematical knowledge in different institutions. The main accent in Ukrainian method of teaching and learning mathematics is given to the development of advanced methodical systems (systems of methods, forms and tools of teaching and learning of mathematics) and evaluation of their effectiveness via pedagogical experiments. Metodika has more practical objectives than didactics. For example, it aims to bring the answer to the following questions: what to teach (content), how to teach (what methods, 10 organizational forms to use, tools), how to evaluate the teaching/learning results, etc.

10. Пиаже Ж. Избранные психологические труды. М., 1994.

11. Maksymenko S.D. Genetic epistemology of J. Piaget / S.D. Maksymenko // Problems of Modern Psychology : Collection of research papers of Kamianets-Podilskyi Ivan Ohienko National University, G.S. Kostiuk Institute of Psychology at the National Academy of Pedagogical Science of Ukraine / scientific editing by S.D. Maksymenko, L.A. Onufriieva. – Issue 32. – KamianetsPodilskyi : Aksioma, 2016. – P. 7–16.

12. Дубасенюк О.А. Підготовка майбутніх учителів до реалізації педагогічної дії. Матеріали Всеукраїнської науково-практичної конференції з міжнародною участю « Теорія і практика підготовки майбутніх учителів до педагогічної дії », 20-21 травня 2011 р., м. Житомир. – Житомир : Вид-во ЖДУ ім. Івана Франка, 2011. – С. 13 – 18.

Fourth, the evidence of central concept, such as 'resource', was called into question: "Understanding resource as something re-sourcing teachers' activity" can be transferred easily into French, or Italian, but is doesn't work in other languages, such as Portuguese:

> 'Recurso' (in Portuguese) is a word composed by the juxtaposition of the prefix «re» and the noun «curso», the first means repetition and the second a path already used, which is the meaning of the Latin recursus (NEGRI, 2007, p. 9). Therefore, recursar (verb in Portuguese) is unusual to give the same meaning of the verb re-source (in English). For that, we used the verb *reabastecer* ou *realimentar* with the idea of 'source again' (extract of the Portuguese Translation issues report, by Katiane Rocha, Cibelle Assis and Sonia Igliori)

It appeared, then, that there is a need to give a sharper definition of the critical concepts of DAD (see Ruthven's comment, § 3).

Fifth, thinking of *adaptation* instead of *translation* leads to establish links between DAD and frames already existing in a given culture, offering new opportunities to rethink the terrain of DAD; for example:

- The links between the instrumental and the documentational genesis through the lens of the mathematics laboratory for the Italian adaptation (Trouche, Gueudet, Pepin, Maschietto, & Panero, 2020);
- The potential links between the concept of instrumentalisation and the Guided Discovery Approach (Goztoniy, 2019), in the case of the Hungarian translation.

Finally, at this stage of the DAD-Multilingual project, these results appear quite promising regarding its aims: deepening the DAD concepts themselves in thinking their possible instantiations in different languages; questioning the translating process itself; questioning the notion of resource itself, resource for/from teaching. These results concern the 'resource' approach to mathematics education, beyond the community of mathematics education, and, beyond, the scientific fields interest in teacher education, cross-cultural studies, and translating processes.

## 5. Perspectives

I would like, as a conclusion, to imagine some perspectives, to be discussed with the actors of the project, perspectives internal to each language, crossing the languages or at a general level.

Each translation-adaptation could live its life in different natural ecosystems: being published in journals, discussed in scientific communities, being used in various research projects (for analyzing teachers interactions with resources, particularly their naming systems, in a variety of contexts), or teaching programs, and crossed with existing approaches. A given language could correspond to diverse cultural, national, or social contexts, e.g. the Spanish adaptation (in Mexico vs. in Spain), the Portuguese adaptation (in Brazil vs. in Portugal), or the Arabic one (in Lebanon vs. in Algeria or Morocco), English constituting a specific case (UK, USA, Australia or India). And, for a given language and a given country, the research context (at University level) and the teaching context (at schools levels) could provide different ecosystems where words and notions may follow their own trajectories. These appropriation processes, at a larger scale, would lead probably to updating each translation, and new issues to be addressed to/by the 'original' DAD frame.

The existing collection of translations opens also different perspectives; for example:

- Confronting the translating techniques, as detailed by Quigley et al. (2011) in the case of the English to Chinese translation:

"*borrowing* (the source language word is transferred directly to the target language), *literal translation* (word-for-word translation), *transposition* (translating the words while paying attention to linguistic differences such as placement of adjectives before or after nouns), *modulation* (a technique often adopted when literal or transposition translation results in a utterance that, though grammatically correct, appears abnormal or awkward), and equivalence (a technique similar to modulation often used in idioms, proverbs, and phrases) […] in order to accurately translate documents, all these techniques must be used. Furthermore, translation through modulation and equivalence requires great attention to cultural, lexical, grammatical, and syntactic aspects of the text.

- Using the diversity of Translating issues reports for a mutual enrichment of each of them, leading towards an updated version of these reports and of the related translations;
- Organizing a new stage in interactions between pairs of languages using their proximity, or origin, for example Spanish-Portuguese, French-Italian, German-Norwegian, Chinese-Japanese, Greek-Ukrainian, Arabic- Hebrew…;
- Using the translated frame for analyzing data in the corresponding language;
- Proposing a special issue of a journal in mathematics education dedicated to the project, its results, issues and perspectives;
- Organizing working groups around concepts, leading towards a common glossary.

Finally, we could imagine a last stage of coming back to the original English presentation of the DAD entry, for updating and perhaps enlarging its scope in the perspective of the 'Resource Approach to Mathematics Education', as described by Gueudet, Pepin, & Trouche (2019).

These perspectives, of course, have to be questioned and enriched by the scientific committee, and discussed within the actors of the project, emerging as a community of (conceptual) enquiry (Jaworski, 2005). In this time of pandemic, sheltering at home, and cultivating a regular social distance, I would like to conclude with a personal statement: such a project, crossing the linguistic and cultural boarders, opens (at least for me) a breathing space, allowing to re-source, really, our conceptualization of teachers' work.

**Acknowledgments…**

- To Ghislaine Gueudet, Birgit Pepin, Ana Isabel Sacristan, and Tommy Dreyfus for their careful rereading and helpful advices;
- To all the participants to the DAD-Multilingual project.